\documentclass{amsart}
\usepackage{amscd}             %Latex2e

\def\C{{\mathbb C}}
\def\H{{\mathbb H}}

\def\R{{\mathbb R}}
\newtheorem{Theorem}{Theorem}
\newtheorem{theorem}{Theorem}[section]
\newtheorem{lemma}[theorem]{Lemma}
\newtheorem{proposition}[theorem]{Proposition}

\newtheorem{conjecture}{Conjecture}
\theoremstyle{definition}

\begin{document}

\title{Asymptotic values of hyperbolic monopoles}
\author{Paul Norbury}
\address{Department of Mathematics and Statistics\\The University of
Melbourne\\Victoria, Australia 3010.}
\email{norbs@ms.unimelb.edu.au}

\begin{abstract}
We show that many hyperbolic monopoles can be distinguished from each
other via their asymptotic values in contrast to the case of Euclidean
monopoles.
\end{abstract} 
\maketitle

\section{Introduction}
\label{sec:intro} 

Magnetic monopoles initially arose out of Dirac's study of the quantum
theory of electro-magnetism.  They are singular solutions of Maxwell's
equations valid away from their singularities.
Bogomolny-Prasad-Sommerfield monopoles are a generalisation of Dirac
monopoles to non-abelian theories where the singularities can be
smoothed away.  They are solutions to the equation $d_A\Phi=*F_A$
where $(A,\Phi)$ is a pair given by a connection with $L^2$ curvature
$F_A$ and a Higgs field---a section of the adjoint bundle---defined on
a trivial bundle over $\R^3$ with structure group a compact Lie group
$G$, and the Hodge star is given with respect to a metric on $\R^3$.
The Higgs field $\Phi$ is constrained to lie in a given orbit of the
Lie algebra on the sphere at infinity.  The mass of the monopole is
defined to be the conjugacy class of the Higgs field on the sphere at
infinity, or more generally it is the orbit in the Lie algebra under
the action of the group of the Higgs field on the sphere at infinity.
The mass also defines an embedding of the circle (the gauge group of
the abelian theory) into the group $G$ so that BPS monopoles can be
compared to Dirac monopoles.

It is reasonable to ask what BPS monopoles look like from a long
distance, and to the extreme, on the sphere at infinity.  It is not
necessarily true that they should look like Dirac monopoles since the
latter approximate BPS monopoles not only far from the singularities
but also only when the singularities are far apart.  We can ask this
question for different metrics on $\R^3$.  When the metric is
Euclidean, the monopoles on the sphere at infinity do look exactly
like Dirac monopoles.  Moreover, as for Dirac monopoles, up to charge,
all Euclidean monopoles look the same at infinity.

In this paper we will consider hyperbolic monopoles defined over
hyperbolic space $\H^3$.  We will show that on the sphere at infinity,
the BPS hyperbolic monopoles take on many different values in contrast
to the Euclidean case.  This agrees with the conjecture that in fact
hyperbolic monopoles are determined by their values on the sphere at
infinity.  This conjecture has been confirmed for integral mass
$SU(2)$ hyperbolic monopoles by Austin and Braam \cite{ABrBou} and it
is an easy fact for hyperbolic Dirac monopoles.  The proof of the
integral mass $SU(2)$ case by Austin and Braam in uses some beautiful
algebraic geometry and introduces discrete Nahm data.  This approach
has been generalised to $SU(n)$ by Murray and Singer \cite{MSiCom} and
it is likely this will lead to a proof that these monopoles are
determined by their asymptotic values.  Our aim is to complement this
work with proofs that work for all gauge groups and non-integral mass.
We are also interested in how such proofs fail for Euclidean
monopoles.

Murray and Singer \cite{MSiSpe} also study the twistor theory of
hyperbolic monopoles.  Their results hold for any mass since rather
than working with circle invariant instantons over the four-sphere
which requires integral mass they work with instantons over Minkowski
space invariant under translations.  They show that a hyperbolic
monopole is determined by its asymptotic value plus some extra
information (remark (3) on p.989.)

It may end up that methods of algebraic geometry will be needed to
show that a hyperbolic monopole is determined by its asymptotic value.
We believe the main result of this paper is still valuable since it
directly shows why hyperbolic monopoles and Euclidean monopoles behave
differently.

As mentioned, a monopole has a mass given by an element of the Lie
algebra, or really the orbit of the element.  We can parametrise the
moduli space of monopoles with given mass by holomorphic maps from the
two-sphere into the orbit in the Lie algebra---see
Section~\ref{sec:hol}.  Since many different masses have isomorphic
orbits and hence the same parametrisation we can speak of a subset of
the moduli space before specifying the mass precisely.  This same idea
is used in the study of ``monopole clouds''.

\begin{Theorem}    \label{th:main}
Given two disjoint compact subsets in the parameter space
of monopoles, if the mass is small enough then the asymptotic values
of the corresponding hyperbolic monopoles respectively give two
distinct subsets.
\end{Theorem}
Remarks: (i) It is probably true that only one of the subsets of the
parameter space need be compact.  We discuss this in
Section~\ref{sec:app}

(ii) One satisfying aspect of the theorem is that we can see where the
proof fails for Euclidean monopoles.

In Section~\ref{sec:dirac} we describe Dirac monopoles over Euclidean
and hyperbolic spaces.  In Section~\ref{sec:hol} we give the
background to the holomorphic map associated to a monopole.  In
Section~\ref{sec:app} we give the proof of Theorem~\ref{th:main}.  We
contrast properties of the asymptotic values of hyperbolic monopoles
with those of Euclidean monopoles in Section~\ref{sec:com}.

\section{Dirac monopoles.}     \label{sec:dirac}
In this section we study the elementary issue of Dirac monopoles since
they give an analogue to the non-linear problem.  It is interesting to
the note that at infinity Euclidean BPS monopoles look exactly like
Euclidean Dirac monopoles whereas hyperbolic space detects a
difference between BPS monopoles and Dirac monopoles.

A Dirac monopole is a solution to Maxwell's equations: $F_A=*d_A\Phi$
where $F_A$ is the magnetic field, or the curvature of a connection
$A$, and $d_A\Phi$ is the electric field given by the covariant
derivative of the Higgs field.  The Hodge star $*$ depends on the
metric. It follows that the magnetic field is harmonic, so for the
Euclidean metric the magnetic field is given by 
\[ B=\frac{\hat{r}}{4\pi r^2}\]
where $\hat{r}$ represents the imaginary-valued 2-form
\[ \hat{r}=\frac{2r^2d\bar{w}dw}{(1+|w|^2)^2},\] 
and more generally 
\[ B=\frac{\widehat{(r-a)}}{4\pi (r-a)^2}\] 
is a monopole at the point $a\in\R^3$.  As $r\rightarrow\infty$,
$r^2B\rightarrow\hat{n}$, the unit normal, which is in particular
independent of $a$.  Thus, on the sphere at infinity a monopole looks
like a symmetric distribution, or the imaginary-valued 2-form,
$2d\bar{w}dw/(1+|w|^2)^2$.  A collection of $k$ monopoles is simply
the sum of $k$ of these and thus gives
\begin{equation}
F_A=\frac{2kd\bar{w}dw}{(1+|w|^2)^2}
\end{equation}
on the sphere at infinity.

Thus we see that (except for the charge) Euclidean Dirac monopoles
cannot be distinguished from a distance.  In contrast with this, the
situation is exactly the opposite for hyperbolic Dirac monopoles.

\begin{proposition}  \label{th:diruniq}
Hyperbolic Dirac monopoles are determined by their asymptotic values.
\end{proposition}
{\em Proof.}  A hyperbolic Dirac monopole is a solution of the
equation $F_A=*d_A\Phi$ for the hyperbolic Hodge star.  The
fundamental solution is 
\[ B=\frac{\hat{r}}{4\pi\sinh^2(r)}\] 
so a more general single Dirac monopole is given by
\[ B_a(x)=\frac{\hat{\nu}}{4\pi\sinh^2(d(x,a))}\] 
where $d(x,a)$ is the hyperbolic distance between a point $x\in\H^3$
and a given point $a$ and $\hat{\nu}$ is the unit vector pointing
(away from $a$) along the geodesic joining $x$ and $a$.  The
asymptotic value of each of these monopoles is given by the unit
outward normal vector of $S^2_{\infty}$ scaled by
$\lim_{r\rightarrow\infty}\sinh^2(r)/\sinh^2(d(x,a))$ (for
$r=d(x,0)$.)  It uniquely determines the monopole since it simply
gives the symmetric measure $d\bar{w}dw/(1+|w|^2)^2$ transformed by
the conformal transformation of $S^2_{\infty}$ induced by the isometry
of $\H^3$ that takes $0$ to $a$.

The most general Dirac monopole is a linear combination of these
single mono\-poles.  The content of this proposition is to show that the
linear combination of conformal transformations of the symmetric
measure on $S^2_{\infty}$ determines the conformal transformations.

A conformal transformation takes 
\[ w\mapsto\frac{aw+b}{cw+d}.\]  
Since the subgroup $SU(2)$ fixes the symmetric measure we need only
consider conformal transformations of the form
\[w\mapsto a_j(w-w_j),\  a_j\in\R^+,\ w_j\in\C.\]  
The symmetric measure maps to $a_j^2d\bar{w}dw/(1+a_j^2|w-w_j|^2)^2$
so a general Dirac monopole has measure at infinity given by
\begin{equation}  \label{eq:meas}
\sum_j \frac{a_j^2d\bar{w}dw}{(1+a_j^2|w-w_j|^2)^2}
\end{equation}
where there might be repeated appearances of a pair $(a_j,w_j)$.

The denominator for (\ref{eq:meas}) is given by
$\Pi_j(1+a_j^2|w-w_j|^2)^2$.  Put $w=x+iy$ and set $y=0$. Then the
denominator factorises with factors $x-x_j\pm i\sqrt(y_j^2+1/a_j^2)$.
Thus, the measure determines each $x_j$ and $y_j^2+1/a_j^2$.  (We have
analytically continued $x$ to take on complex values.)  Similarly, if
we set $x=0$ then we get each $y_j$ and $x_j^2+1/a_j^2$ so we get each
$x_j,y_j,a_j$ and the boundary measure has determined the Dirac
monopole.\qed

The space of compactly supported continuous functions on hyperbolic
space acts on the symmetric measure $\omega$ on $S^2$ by
$f\cdot\omega=\int_{\H^3}f(x)x\cdot\omega dx$ where we think of $x\in
SL(2,\C)/SU(2)$.  Atiyah has suggested that this action might be faithful.

\begin{conjecture}
Let $f\in C^0_c(\H^3)$, then $f\cdot\omega=0\ \Leftrightarrow\ f\equiv 0.$
\end{conjecture}
The proposition would then fit into a rather natural setting,
following from limiting behaviour of such a result.

\section{Holomorphic maps.}    \label{sec:hol}
In this section we will describe the holomorphic map of the two-sphere
into a homogeneous space associated to a monopole via scattering.
First we will give a brief description of the homogeneous spaces.

Let $\xi\in{\bf g}$, the Lie algebra of $G$, and let $K_{\xi}=\{ g\in
G\ |\ g\cdot\xi=\xi\}$ be the isotropy subgroup of $\xi$ from the
adjoint action of $G$.  Then the homogeneous space $G/K_{\xi}$ is a
complex manifold and we can speak of holomorphic maps into this
manifold.  This is best seen using the isomorphism $G/K_{\xi}\cong
G\cdot\xi\subset{\bf g}$.

The tangent space of an adjoint orbit $X=G\cdot\xi\subset{\bf g}$ has
a nice description.  At $\eta=g\cdot\xi\in G\cdot\xi$, the tangent
space $T_{\eta}X=[{\bf g},\eta]\subset{\bf g}$.  It is more convenient
to use ${\bf g}/\ker [\eta,\cdot]\cong T_{\eta}X$.  The isomorphism is
given by $u\mapsto [u,\eta]$ for $u\in{\bf g}$.  

Homogeneous spaces have complex realisations 
\begin{equation}  \label{eq:hom}
G/K_{\xi}\cong G^c/P_{\xi}
\end{equation} 
where $P_{\xi}$ is the parabolic subgroup of $G^c$ with the further
property that $P_{\xi}\cap G=K_{\xi}$.  The isomorphism (\ref{eq:hom})
simply says that given any $g\in G^c$, there exists $p\in P_{\xi}$
such that $gp\in G$ and $p$ is unique up to $p\mapsto pk$ for $k\in
K_{\xi}$.  When $G=U(n)$, this is the Gram-Schmidt process.

The complex structure at $\xi$ is given by $Ju\equiv iu({\rm mod\ }{\bf
p}_{\xi})$ with respect to the trivialisation $T_{\xi}X\cong{\bf
g}/\ker [\xi,\cdot]$.  It is well-defined since given $iu$, there is
an element $v$ of ${\bf p}_{\xi}$, the Lie algebra of $P_{\xi}$,
unique up to an element of $\ker [\xi,\cdot]$ such that $iu+v\in{\bf
g}$.  The complex structure at each point of the orbit $\eta=g\cdot\xi$
is defined similarly.

A map $f:S^2\rightarrow G/K_{\xi}$ is holomorphic when its lift $u$ to
$G$ (defined locally) satisfies
$u^{-1}\partial_xu+Ju^{-1}\partial_yu=0$, or equivalently for $w=x+iy$
\begin{equation}   \label{eq:holo}
u^{-1}\partial_{\bar{w}}u(w)\subset{\bf p}_{\xi}.
\end{equation}

The adjoint orbit has a natural symplectic structure, compatible with
the complex structure to give a Kahler structure, given at $T_{\eta}X$
by $\omega(u,v)=\langle\eta,[u,v]\rangle$ where
$\langle\cdot,\cdot\rangle$ is the Killing form.  At $T_{\eta}X$, the
metric is $g(u,v)=\omega(Ju,v)$ for the complex structure $J$.  There
are many other symplectic structures and complex structures that arise
less naturally.

The map $f:S^2\rightarrow G/K_{\xi}$ pulls back the symplectic form
$\omega$ to a two-form over $S^2$ via its lift $u$:
\begin{equation}
f^*\omega=\langle\xi,[u^{-1}\partial_wu,u^{-1}\partial_{\bar{w}}u]\rangle 
d\bar{w}dw.
\end{equation}

A hyperbolic monopole $(A,\Phi)$ with finite energy,
$\|F_A\|_2<\infty$, has a well-defined limit at infinity (ensuring
that the problem of this paper is well-posed) and the components of
the monopole and their derivatives satisfy asymptotic decay conditions
near infinity, \cite{MRaPri}.  In particular,
$\Phi|_{S^2_{\infty}}=\xi\in{\bf g}$ is the mass of the monopole.

More precisely, there exists a gauge in which 
\[ \Phi=\xi+O(e^{-cr})\]
and $A=A_wdw+A_{\bar{w}}d\bar{w}+A_r$ such that
\[ A_wdw+A_{\bar{w}}d\bar{w}=A_{\infty}+O(e^{-cr})\] 
for a $K_{\xi}$ connection $A_{\infty}$ on $S^2$, and 
\[A_r=O(e^{-cr})\] 
as $r\rightarrow\infty$ and $c>0$ is a constant.  There are similar
estimates on the derivatives.

The asymptotic conditions on the monopole ensure along each radial
geodesic the existence of a frame of fundamental solutions
$g:\R^+\rightarrow G^c$ of the scattering equation
\begin{equation}    \label{eq:scat}
(\partial_r^A-i\Phi)g=0
\end{equation}
with the property that $g(0)\in G$ and
$\lim_{r\rightarrow\infty}g\exp(-i\xi r)$ is bounded.  The solution is
unique up to $g\mapsto gk$ for $k\in K_{\xi}$.  

We can choose a family of solutions $g(w,r)$ to (\ref{eq:scat}) that
depend smoothly on $w$.  It follows from the Bogomolny equation
$F_A=*d_A\Phi$ that
\[ (\partial_r^A-i\Phi)\partial_{\bar{w}}^Ag=0\] 
and $\partial_{\bar{w}}^Ag\exp(-i\xi r)$ is bounded as
$r\rightarrow\infty$, so 
\begin{equation}  \label{eq:gauge}
\partial_{\bar{w}}^Ag=g\eta(w)
\end{equation}
for some $\eta(w)\in{\bf p}_{\xi}$.  In particular,
$\partial_{\bar{w}}g(w,0)=g(w,0)\eta(w)$ since by the choice of
coordinate system $\partial_{\bar{w}}^A=\partial_{\bar{w}}$ at $r=0$.
Thus $g(w,0)^{-1}\partial_{\bar{w}}g(w,0)\subset{\bf p}_{\xi}$ and by
(\ref{eq:holo}) this means $g(w,0):S^2\rightarrow G/K_{\xi}$ is a
holomorphic map.  In order to make sense of the value of $g$ at $r=0$
we have chosen a frame of the bundle there.  This construction gives
part of the following theorem.
\begin{theorem}\cite{JNoDeg}  \label{th:diff}
The space of hyperbolic monopoles framed at $0\in\H^3$ with gauge
group $G$ and mass $\xi$ is diffeomorphic to the space of holomorphic
maps $Hol(S^2,G/H)$ where $H$ is the isotropy subgroup of $\xi$.
\end{theorem}

We can interpret the solution $g$ of (\ref{eq:scat}) as a choice of
gauge and then (\ref{eq:scat}) and (\ref{eq:gauge}) give $(A,\Phi)$
with respect to this gauge, respectively showing that $A_r-i\Phi=0$
and $A_{\bar{w}}=\eta(w)$.  We can choose another solution
$g(w,r)p(w)$ of (\ref{eq:scat}) for $p:\C\rightarrow P_{\xi}$ that has
the same asymptotic properties as $g$ but no longer satisfies
$gp(w,0)\in G$ and with the property that with respect to this gauge
\begin{equation}  \label{eq:van} 
A_r-i\Phi=0,\ A_{\bar{w}}=0.
\end{equation}
This simply uses the fact that any holomorphic map $f:S^2\rightarrow
G/K_{\xi}$ (locally) lifts to a map $u:\C\rightarrow G$, and there is
a map $p:\C\rightarrow P_{\xi}$ such that $up:\C\rightarrow G^c$ is a
lift of $f$ to an algebraic map.  The maps $u$ and $p$ are unique up
to $(u,p)\mapsto (uk,k^{-1}p)$ for $k:\C\rightarrow K_{\xi}$.  The
evaluation $g(w,0)$ is of course a lift of a holomorphic map to
$G/K_{xi}$.

Since $(A,\Phi)$ is Hermitian, with respect to the frame $gp$
satisfying (\ref{eq:scat}), the Hermitian metric
$H=(gp)^*gp:\H^3\rightarrow G^c/G$ together with (\ref{eq:van}) gives
the remainder of the monopole
\begin{equation} \label{eq:van2}
A_r+i\Phi=H^{-1}\partial_rH,\ A_w=H^{-1}\partial_wH.
\end{equation}
The Bogomolny equations become $B(H)=0$ where
\begin{equation}  \label{eq:hermbog}
B(H)=\sinh^2(r)\partial_r(H^{-1}\partial_rH)+(1+|w|^2)^2\partial_{\bar{w}}
(H^{-1}\partial_wH).
\end{equation}
See \cite{JarMon} for further details.  

Notice that $H(0)=p(w)^*p(w)$ is not well-defined (it depends on $w$)
and $H|_{S^2_{\infty}}$ is reduced i.e. $H|_{S^2_{\infty}}\in
K_{\xi}^c/K_{\xi}$ (for generic mass $K_{\xi}$ is a torus and
$H|_{S^2_{\infty}}$ is a potential.)

Now we represent a monopole as a Hermitian metric $H$ that satisfies
(\ref{eq:hermbog}).  Given $H_1$ and $H_2$, define the endomorphism
$h=H_1^{-1}H_2$.
\begin{lemma}
If $\Phi_1=\Phi_2$ on $S^2_{\infty}$, so the two monopoles have the
same mass and we choose gauges in which the Higgs fields look the
same, then the endomorphism $h$ is conjugate to a bounded endomorphism.
\end{lemma}
{\em Proof.}  The Hermitian metric $H_j$ arises from $(A_j,\Phi_j)$ as
$H_j=g_j^*g_j$.  Put $\Phi_j|_{S^2_{\infty}}=\xi$, then
$g_j=G_j(w,r)exp(i\xi r)$ for $G_j(w,r)$ bounded so $g_2g_1^{-1}$ is
bounded.  Now, $h=g^{-1}(((g_2g_1^{-1})^*g_2g_1^{-1})g_1$.
\qed\\

The two monopoles have the same asymptotic value precisely when
$h|_{S^2_{\infty}}=I$, the identity endomorphism.  The complete metric
on the space of Hermitian metrics given by taking the supremum over
hyperbolic space of $\langle H^{-1}\delta H,H^{-1}\delta H\rangle$
uses the Killing form so the previous lemma implies that for two
monopoles with the same mass, the distance between $H_1$ and $H_2$ is
finite.  Two monopoles are the same when the distance between their
Hermitian metrics is zero.

\section{Approximate monopoles.}       \label{sec:app}
In this section we will prove Theorem~\ref{th:main}.  Our strategy is
as follows.  For each holomorphic map and mass we can find an
approximate monopole and a unique exact monopole nearby.  The smaller
the mass, the better the approximation.  For any two holomorphic maps,
the distance between the asymptotic values of the corresponding
approximate monopoles is independent of the mass and positive.  Thus,
for small enough mass, when the two approximations are quite good,
there must also be a positive distance between the asymptotic values
of the two exact monopoles and the theorem is proven.  We will
actually use Hermitian metrics in place of monopoles since there is a
good notion of distance between Hermitian metrics and there are
techniques to estimate this distance.

A Hermitian metric $H:\H^3\rightarrow G^c/G$ can be associated to a
more general set of pairs $(A,\Phi)$ than monopoles over $\H^3$.  In
fact, to any pair $(A,\Phi)$ that satisfies
\[[\partial_{\bar{w}}^A,\partial_r^A-i\Phi]=0\] 
we can associate a Hermitian metric $H$ and $(A,\Phi)$ is retrieved
from $H$ by (\ref{eq:van}) and (\ref{eq:van2}).  This is the class of
pairs we will consider.

Given $\xi\in\bf g$ and a holomorphic map $f:S^2\rightarrow G/K_{\xi}$
define 
\begin{equation}   \label{eq:appherm}
{\mathcal H}_f=p^*\exp(2i\xi r)p
\end{equation} 
where $p:\C\rightarrow P_{\xi}$ is a map into the parabolic subgroup
with the property that $up$ is a lift of $f$ to an algebraic map from
$\C$ to $G^c$ and $u$ is a lift of $f$ to a map from $\C$ to $G$.  The
map $u$ is ambiguous up to an action of $K_{\xi}$ on the right and
hence $p$ inherits this ambiguity on the left.  The expression for
${\mathcal H}_f$ is independent of this ambiguity.

In order to show the existence of a monopole for any given mass and
holomorphic map we use the non-linear heat flow for Hermitian metrics
with initial value given by ${\mathcal H}_f$.
\begin{equation}   \label{eq:flow} 
H^{-1}(w,r,t)\partial H(w,r,t)/\partial t=B(H(w,r,t)),\ 
H(w,r,0)={\mathcal H}_f(w,r)
\end{equation}
\begin{theorem}\cite{JarMon,JNoDeg}
There is a unique solution $H(w,r,t)$ of (\ref{eq:flow}).
\end{theorem}
The solution $H(w,r,t)$ of the heat flow converges to a Hermitian
metric that satisfies $B(H(w,r,\infty))=0$ and gives rise to a
monopole with holomorphic map $f$.  Together with the scattering
construction described in Section~\ref{sec:hol} this gives the proof
of Theorem~\ref{th:diff}.

It is worth pointing out that the construction of monopoles from
holomorphic maps is treated differently in \cite{JarMon} and
\cite{JNoDeg} and here it is treated slightly differently again.
In \cite{JarMon} the initial choice of Hermitian metric used
explicitly known symmetric hyperbolic monopoles.  In \cite{JNoDeg},
since both instantons and hyperbolic monopoles were treated together
it was more convenient to choose an initial Hermitian metric that was
independent of such information (and also to use something more
general than a Hermitian metric.)  Neither of these suffice for our
purposes here.  In order that the limiting connections at infinity of
different monopoles can be compared we need to ensure that a common
reduction of the monopoles to a subgroup (usually a maximal torus) is
used.  This is why the parabolic subgroup is specified and features in
the Hermitian metric above.  In particular, the approximate monopole
defined by ${\mathcal H}_f$ has the same asymptotic mass as the
monopole associated to $H(w,r,\infty)$, so $d(H(w,r,\infty),{\mathcal
H}_f)$ is finite.

The metric on the space of Hermitian metrics is given by
$tr(H^{-1}dH)$ so the heat flow gives an estimate of the distance from
the initial ${\mathcal H}_f$ and the final $H(w,r,\infty)$:
\begin{eqnarray}
d(H(w,r,\infty),{\mathcal H}_f)&\leq&\int_0^{\infty}|B(H(w,r,t))|dt\\
&\leq&\int_0^{\infty}{\rm max}_{|\rho|=s}|B(\mathcal H_f(w,\rho))|^2G(s,r)ds
\label{eq:est}
\end{eqnarray}
where the second inequality comes from the fact that $|B(H(w,r,t))|$
is dominated by a solution of the linear heat flow for a Laplacian
like operator that reduces to the usual Laplacian on radially
symmetric functions.  When we maximise $|B(\mathcal H_f)|$ over
spheres of constant radius we get a function of the radius so we can
use a one-dimensional Green's function $G(s,r)$.  See
\cite{JarMon,JNoDeg} for the proof of this and also
\cite{DonAnt,SimCon} where this technique is introduced.

The following two propositions estimate how well (\ref{eq:appherm})
approximates a monopole by using (\ref{eq:est}).  The first
proposition is enough to prove Theorem~\ref{th:main}.  We go on to
prove more in the second proposition.  It relates the estimate in the
first proposition with $|f^*\omega|$, the two-form on $S^2$ pulled
back by the holomorphic map, where $f^*\omega$ is compared to the
standard two-form on $S^2$ to get its magnitude.  This precise
information is included particularly to show how good the
approximation is as the holomorphic maps bubble.

\begin{proposition}   \label{th:dist}
$d(H(w,r,\infty),{\mathcal H}_f)\leq C|\xi|\ {\rm max}_{S^2_{\infty}}|F_A|$.
\end{proposition}
\begin{proof}
By (\ref{eq:est}) we have to get an upper bound for
$|B({\mathcal H}(w,r))|$ on each sphere $r=$ constant.  Using
(\ref{eq:van}) and (\ref{eq:van2}) with $H={\mathcal H}_f$, and going
to a unitary gauge we get
\[ \Phi=\xi, A_r=0, A_{\bar{w}}=\exp(i\xi r)u^{-1}\partial_{\bar{w}}u
\exp(-i\xi r), A_w=\exp(-i\xi r)u^{-1}\partial_wu\exp(i\xi r)\]
and $B({\mathcal H}(w,r))=-i(1+|w|^2)^2F_{\bar{w}w}$.

The connection $A$ splits into a radially independent $K_{\xi}$
connection and an exponentially decaying connection.  More precisely,
put
\[{\bf g}^c={\bf k}_{\xi}\oplus{\bf n}^+_{\xi}\oplus{\bf n}^-_{\xi}\]
where ${\bf k}_{\xi}=\ker [\xi,\cdot]$ is the Lie algebra of $K_{\xi}$
and ${\bf k}_{\xi}\oplus{\bf n}^+_{\xi}={\bf p}_{\xi}$, the Lie
algebra of $P_{\xi}$.  Alternatively, we can characterise the sub-spaces by
\[(\exp{(i\xi r)}-I)\cdot{\bf k}_{\xi}=0=\lim_{r\rightarrow+\infty}
\exp{(\pm i\xi r)}\cdot{\bf n}^{\pm}.\]
Decompose $v\in{\bf g}^c$ accordingly as $v=v^0+v^++v^-$.

The connection $A$ decomposes as $A=A^0+a$ for
$a=A_{\bar{w}}^+d\bar{w}+A_w^-dw$ with the property that $A^0$ is a
$K_{\xi}$ connection independent of $r$, and $a$ is a 1-form that
decays exponentially as $r\rightarrow\infty$. Then
$F_{\bar{w}w}=F_{A^0}+d_{A^0}a+a\wedge a$ and $F_{A^0}$ is independent
of $r$ whilst the rest decays exponentially, so
$F_{A^0}=F_A|_{S^2_{\infty}}$.

Define $c_{\xi}>0$ to be the smallest eigenvalue of the action of
$\exp(i\xi r)$ on ${\bf g}^c$.  Then each time we say that a term
decays exponentially, it decays at least as fast as $e^{-c_{\xi}r}$.
Notice that $c_{\xi}\leq|\xi|$.  We have $|F_A|\leq
M_1|F_{A^0}|(1-e^{-c_{\xi}r})$ for some constant $M_1\geq 1$ and thus
\begin{eqnarray*}
d(H(w,r,\infty),{\mathcal H}_f)&\leq&\int{\rm max}_{|\rho|=s}|B(\mathcal
H)(w,\rho)|^2G(s,r)ds
\\ &\leq&M_1{\rm max}_{S^2_{\infty}}|F_A|\int_0^{\infty}\frac{(1-e^{-c_{\xi}s})
{\rm min}\{ r,s\}}{\sinh^2s}ds
\end{eqnarray*}
where ${\rm min}\{ r,s\}/\sinh^2s$ is the one-dimensional Green's
function.

Now 
\[\int_0^{\infty}\frac{(1-e^{-c_{\xi}s}){\rm min}\{ r,s\}}{\sinh^2s}ds=
\int_0^r\frac{(1-e^{-c_{\xi}s})s}{\sinh^2s}ds+r\int_r^{\infty}
\frac{(1-e^{-c_{\xi}s})}{\sinh^2s}ds\]
and the second term of the right hand side converges to $0$ as
$r\rightarrow\infty$.  Since $1-e^{-c_{\xi}s}\leq c_{\xi}s$, the first
term is dominated by $c_{\xi}M_2$ for a constant $M_2$.  Since
$c_{\xi}\leq|\xi|$ the proposition follows.  The constant $C$ in the
statement of the proposition does depend on the holomorphic map $f$,
and is bounded below by a constant independent of $f$.
\end{proof}

In particular, the estimate depends only on the holomorphic map and
the mass.  For small mass, the distance is small.  For any two
holomorphic maps $f$ and $g$, notice that restricted to $S^2_{\infty}$
the distance $d({\mathcal H}_f,{\mathcal H}_g)$ is independent of the
mass since it depends on
\[{\mathcal H}_f^{-1}{\mathcal H}_g=p_f^{-1}\exp(-2i\xi r)(p_f^*)^{-1}
p_g^*\exp(2i\xi r)p_g\] which is independent of $\xi$ in the limit
$r\rightarrow\infty$.  Of course the holomorphic maps $f$ and $g$ use
$\xi$ to be defined, but they could equally well use a mass that gives
the same homogeneous manifold, like $\lambda\xi$ for $\lambda\in\R^*$.

Thus, if we take any two holomorphic maps and choose the mass small enough,
then the nearby monopoles must have different asymptotic values.  This
argument extends to two disjoint compact subsets of the space of
holomorphic maps and Theorem~\ref{th:main} is proven.

\begin{proposition}
$\langle F_A,\Phi\rangle|_{S^2_{\infty}}=f^*\omega.$
\end{proposition}
\begin{proof}
In the proof of the previous proposition we saw that $F_A$ has a
radially independent part $F_A|_{S^2_{\infty}}$ and an exponentially
decaying part.  Since $F_{\bar{w}w}$ vanishes at $r=0$ we get an
identity relating $F_A|_{S^2_{\infty}}$ with the term that cancels it.

In the notation of the previous proof, 
$F_{\bar{w}w}^0d\bar{w}dw=F_A|_{S^2_{\infty}}+[A_{\bar{w}}^+,A_w^-]^0
d\bar{w}dw$ so 
\[\langle F_A,\Phi\rangle|_{S^2_{\infty}}=\langle -[A_{\bar{w}}^+,A_w^-]^0
d\bar{w}dw,\Phi\rangle|_{r=0}=\langle -[A_{\bar{w}}^+,A_w^-]d\bar{w}dw,
\Phi\rangle|_{r=0}\]
where the last equality follows from the fact that $\Phi$ is
orthogonal to ${\bf n}^{\pm}$.  But at $r=0$,
$A_{\bar{w}}=u^{-1}\partial_{\bar{w}}u, A_w=u^{-1}\partial_wu$ so
\[\langle F_A,\Phi\rangle|_{S^2_{\infty}}=-\langle [u^{-1}\partial_{\bar{w}}u,
u^{-1}\partial_wu],\xi\rangle d\bar{w}dw=f^*\omega.\]
\end{proof}

The previous proposition shows that the heat flow gives bad estimates
for well-separated monopoles.  That is, if a sequence of holomorphic
maps bubble then the pull-back of the Kahler form will bubble and
$\langle F_A,F_A\rangle$ which determines the accuracy of the
approximate monopole, gets a contribution from $\langle
F_A,\Phi\rangle$, and hence gets large.  (It may even be true that $\langle
F_A,\Phi\rangle$ controls $\langle F_A,F_A\rangle$.)

Well-separated monopoles are Dirac-like and are the source of our
intuition that hyperbolic monopoles have interesting asymptotic
limits.  So far we have not been able to produce good approximate
well-separated monopoles.  It would be very interesting to see such
approximations since they would combine the linear nature of Dirac
monopoles with the soliton nature of gluing together holomorphic maps.

Given the intuition that asymptotic values of well-separated monopoles
look like asymptotic values of Dirac hyperbolic monopoles, we would be
able to relax the condition in Theorem~\ref{th:main} allowing one set
to be non-compact and in particular apply the theorem to a point and a
deleted neighbourhood of the point.  

\begin{conjecture}
Given a compact subset and a disjoint subset in the parameter space
of monopoles, if the mass is small enough then the asymptotic values
of the corresponding hyperbolic monopoles respectively give two
distinct subsets.
\end{conjecture}

A related and interesting issue is to know if the set of monopoles
with bounded curvature on $S^2_{\infty}$ and bounded mass, gives rise
to a compact set in the space of holomorphic maps.  Such a result
would also prove the conjecture.  In special cases it can be shown
that for a fixed holomorphic map the maximum value of the curvature at
infinity is a monotone decreasing function of the mass so the
conjecture follows.

\section{Comparison with Euclidean monopoles.}  \label{sec:com}

The asymptotic value of a Euclidean monopole gives a symmetric
connection on the sphere at infinity, and in particular, all monopoles
(with the same mass and charge) give rise to the same connection at
infinity.  This contrasts with the hyperbolic case and it is
interesting to see where the point of departure from the behaviour of
hyperbolic monopoles occurs.

The proof of Proposition~\ref{th:dist} goes through for Euclidean
monopoles with the only change being in the one-dimensional Green's
function.  We use ${\rm min}\{ r,s\}/s^2$ instead of ${\rm min}\{
r,s\}/\sinh^2s$.  Now
\[ d(H(w,r,\infty),{\mathcal H}_f)\leq C\ {\rm max}_{S^2_{\infty}}|F_A|
\left(\int_0^r\frac{(1-e^{-c_{\xi}s})}{s}ds+r\int_r^{\infty}
\frac{(1-e^{-c_{\xi}s})}{s^2}ds\right)\]
and the second term is bounded whilst the first term is $O(\ln r)$.
This is enough to show that the heat flow converges and thus
Theorem~\ref{th:diff} is true for Euclidean monopoles,
\cite{JarMon,NorPer}.  But we see that the asymptotic value of the 
monopole can move arbitrarily far and a posteriori we know that all of
the asymptotic values converge to the same connection.

This is an appropriate place to mention the result of Murray and
Singer \cite{MSiSpe} regarding asymptotic values of hyperbolic
monopoles.  They show that an $SU(n)$ hyperbolic monopole is
determined by $(\nabla^0|_{S^2_{\infty}})^{0,1}$ and
$(b_+|_{S^2_{\infty}})^{0,1}$, using their notation.  The first term
is the asymptotic value of the monopole and the second term is an
artifact of the holomorphic gauge they use giving off-diagonal terms.
In the good unitary gauge defined by Rade
\cite{RadSin} the asymptotic value is a $U(1)^{n}$ connection.
The term $(b_+|_{S^2_{\infty}})^{0,1}$ essentially encodes the
holomorphic map which is also enough to give a Euclidean monopole so
no new behaviour is seen there.  It is not so surprising since their
methods are similar to those applied to Euclidean monopoles.

Finally, we mention a maximum principle which a priori may have led to
a proof that hyperbolic monopoles are determined by their asymptotic
values.  It ends up that the maximum principle also applies to
Euclidean monopoles so it proves a result that is true for both cases.

As in Section~\ref{sec:hol} define $h=H_1^{-1}H_2$ for two Hermitian
metrics satisfying $B(H_i)=0$.  Consider
$\sigma(h)=tr(h)+tr(h^{-1})-2n, (n=tr I)$.  This is a non-negative
quantity that vanishes precisely when $H_1=H_2$.  The problem of
showing that a monopole is determined by its asymptotic values then
becomes the problem of showing that if the asymptotic value of
$\sigma(h)$ vanishes then $\sigma(h)$ vanishes identically.  The
following inequality leads to a maximum principle.
\[\sinh^2(r)\partial^2_r\sigma+(1+|w|^2)^2\partial_{\bar{w}}\partial_w\sigma
\geq 0.\]
It applies to Euclidean monopoles also
\[r^2\partial^2_r\sigma+(1+|w|^2)^2\partial_{\bar{w}}\partial_w\sigma\geq 0.\]
The maximum principle states that $\sigma(h)$ is dominated by any
function that dominates $\sigma(h)$ on the boundary and lies in the
kernel of the second order partial differential operator above.  The
important point here is that there are two boundary components, $r=0$
and $r=\infty$, since $\sigma(h)$ depends on $w$ at $r=0$.  The
function $a+br$ is a good comparison function for constants $a$ and
$b$ chosen so that $a\geq{\rm max}_{r=0}\sigma(h)$ and $b>0$.  As
$b\rightarrow 0$ we see that
\[{\rm max}_{r=0}\sigma(h)\geq{\rm max}_{r=\infty}\sigma(h).\]  
This is true of both hyperbolic monopoles and Euclidean monopoles.  In
the latter case, ${\rm max}_{r=\infty}\ \sigma(h)=0$ so the inequality
is trivial.

\end{document}